\documentclass[10pt]{article}

\usepackage{amsmath,amsfonts,amsthm,amssymb,graphicx,qtimes,eucal}

\theoremstyle{plain}
\newtheorem{thm}{Theorem}[section]
\newtheorem{lem}[thm]{Lemma}
\newtheorem{prop}[thm]{Proposition}
\newtheorem{cor}[thm]{Corollary}

\theoremstyle{definition}

\newtheorem{ex}[thm]{Example}

\newtheorem*{qtn}{Question}
\newtheorem*{rem}{Remark}

\theoremstyle{remark}

\newcommand{\lra}{\longrightarrow}
\newcommand{\sig}{\sigma}
\newcommand{\al}{\alpha}
\newcommand{\Z}{\mathbb Z}
\newcommand{\R}{\mathbb R}
\newcommand{\Q}{\mathbb Q}
\newcommand{\A}{\mathbb A}
\newcommand{\gam}{\gamma}
\newcommand{\Sig}{\Sigma}
\newcommand{\lam}{\lambda}

\newcommand{\C}{\mathbb C}
\newcommand{\Hy}{\mathbf H}
\newcommand{\thet}{\theta}
\newcommand{\Gam}{\Gamma}
\newcommand{\conj}{\overline}
\newcommand{\om}{\omega}
\newcommand{\Del}{\Delta}
\newcommand{\mc}{\mathcal}
\newcommand{\mf}{\mathfrak}
\newcommand{\ep}{\epsilon}
\DeclareMathOperator{\M}{M}
\DeclareMathOperator{\SL}{SL}
\DeclareMathOperator{\PSL}{PSL}
\DeclareMathOperator{\GL}{GL}
\DeclareMathOperator{\SO}{SO}
\DeclareMathOperator{\PU}{PU}
\DeclareMathOperator{\SU}{SU}
\DeclareMathOperator{\Sp}{Sp}
\DeclareMathOperator{\PSO}{PSO}
\DeclareMathOperator{\PGL}{PGL}
\DeclareMathOperator{\rank}{rank}
\DeclareMathOperator{\Gal}{Gal}
\DeclareMathOperator{\Id}{Id}
\DeclareMathOperator{\End}{End}
\DeclareMathOperator{\im}{Im}
\DeclareMathOperator{\real}{Re}
\DeclareMathOperator{\Hom}{Hom}

\newenvironment{mat}{\left(\begin{matrix}}{\end{matrix}\right)}

\begin{document}

\title{Property (FA) and lattices in $\SU(2,1)$}
\author{Matthew Stover \\ \small{University of Texas at Austin}\\ \small{\textsf{mstover@math.utexas.edu}}}
\date{\today}

\maketitle

\begin{abstract}
In this paper we consider Property (FA) for lattices in $\SU(2,1)$. First, we prove that $\SU(2,1;\mc{O}_3)$ has Property (FA). We then prove that the arithmetic lattices in $\SU(2,1)$ of second type arising from congruence subgroups studied by Rapoport--Zink and Rogawski cannot split as a nontrivial free product with amalgamation; one such example is Mumford's fake projective plane. In fact, we prove that the fundamental group of any fake projective plane has Property (FA).
\end{abstract}

\section{Introduction}\label{intro}

Two important questions in the study of lattices in semisimple Lie groups are whether a given lattice splits as a nontrivial free product with amalgamation or admits a homomorphism onto $\Z$. Property (FA), originally due to Bass and Serre, encodes precisely when a finitely generated group has neither of these properties (see Theorem \ref{bigfa}). More generally, one can ask for these properties to hold in a finite sheeted cover. The virtual-$b_1$ conjecture asks, most famously in the setting of closed hyperbolic 3-manifolds, whether the fundamental group of a given manifold has a finite index subgroup admitting a homomorphism onto $\Z$. Using Kazhdan's Property (T), one can prove that irreducible lattices in $\Sp(n,1)$ for $n\geq 2$, $F_{4(-20)}$, and semisimple Lie groups with $\R$-rank at least 2 always have Property (FA) (see \cite{HV}). Therefore, all of the interesting questions relating Property (FA) and irreducible lattices in semisimple Lie groups occur for the fundamental groups of real and complex hyperbolic manifolds -- lattices in $\PSO_0(n,1)$ and $\PU(n,1)$.

Due to exceptional isomorphisms, we can consider the fundamental groups of real hyperbolic 2-manifolds and complex hyperbolic 1-manifolds, lattices in $\PSO_0(2,1)$ and $\PU(1,1)$, as Fuchsian groups, i.e.\ lattices in $\PSL(2;\R)$. Splittings as a free product with amalgamation for cocompact Fuchsian groups are well understood; see \cite{MLR} for a complete list of the known results. For example, considering a separating curve on a compact Riemann surface it follows that many finite covolume Fuchsian groups split as nontrivial amalgamated products. Cocompact Fuchsian triangle groups are well known to have Property (FA), but $\PSL(2;\Z)$ -- the $(2,3,\infty)$ triangle group -- is a free product of two finite cyclic groups. Furthermore, it is known that all Fuchsian groups virtually surject onto $\Z$.

However, the situation becomes much more complicated for fundamental groups of real hyperbolic 3-manifolds and orbifolds, considered (again by an exceptional isomorphism) as lattices in $\PSL(2;\C)$. If $d$ is a square free natural number, Frohman and Fine \cite{ff} prove that the Bianchi group $\PSL(2;\mc{O}_d)$ splits as a nontrivial free product with amalgamation for $d\neq 3$, where $\mc{O}_d$ denotes the ring of integers in $\Q(\sqrt{-d})$, and Serre proves in \cite{Se} that $\PSL(2;\mc{O}_3)$ has Property (FA). Using similar techniques to Serre, we prove the following complex hyperbolic analogue.

\begin{thm}\label{zeta3}
$\SU(2,1;\mc{O}_3)$ and $\PU(2,1;\mc{O}_3)$ have Property (FA).
\end{thm}

The relative similarity of the proofs for $\PSL(2;\mc{O}_3)$ and $\PU(2,1;\mc{O}_3)$ begs the question as to how much further this analogy between $\PSL(2;\mc{O}_d)$ to $\PU(2,1;\mc{O}_d)$ carries. A theorem of Shimura \cite{Shi} implies that $\SU(2,1;\mc{O}_d)$ virtually surjects onto $\Z$ for all $d$, though no explicit homomorphisms are known for $d\neq 3$, so we pose:

\begin{qtn}
Does $\SU(2,1;\mc{O}_d)$ or $\PU(2,1;\mc{O}_d)$ have Property (FA) for $d\neq 3$?
\end{qtn}

Theorem \ref{zeta3} indicates that there is a connection between certain real and complex hyperbolic lattices. In fact, complex hyperbolic lattices seem to bridge the gap between the nonrigidity found in real hyperbolic lattices and rigid higher rank phenomenon. Arithmeticity and superrigidity were shown to hold for all irreducible higher rank lattices by Margulis (see \cite[Chap. 0]{Mar}), and a combination of work of Corlette \cite{Cor} and Gromov--Schoen \cite{GS} implies that these properties also hold for irreducible quaternionic hyperbolic lattices. However, superrigidity fails dramatically for real hyperbolic lattices, and non-arithmetic lattices exist in $\SO(n,1)$ for all $n$ \cite{GPS}. In the complex hyperbolic setting, non-arithmetic lattices are known to exist in $\SU(n,1)$ for $n=2,3$ \cite{DM} and whether non-arithmetic lattices exist in $\SU(n,1)$ for $n\geq 4$ remains a major open question.

When arising from congruence subgroups, arithmetic lattices in $\SU(2,1)$ of the second type have several properties that are remarkably similar to those of superrigid lattices. These include non-archimedean and archimedean superrididity-like properties and vanishing first cohomology -- see Sec. \ref{secondproof}. As Rogawski proves in \cite{Ro}, these lattices have $b_1=0$, and Blasius and Rogawski prove in \cite{BR} that these lattices have Picard number one. In fact, it is a question attributed to Rogawski as to whether all lattices in $\SU(2,1)$ satisfying these criteria are necessarily arithmetic and of the second type. We strengthen the superrigid-like analogy for these lattices with the following theorem, which is the primary result of this paper.

\begin{thm}\label{picard}
If $\Gam<\SU(2,1)$ is a torsion-free congruence arithmetic lattice of second type, then $\Gam$ does not split as a nontrivial free product with amalgamation.
\end{thm}

Since these lattices have $b_1=0$, this immediately implies that $\Gam$ has Property (FA). In particular, Theorem \ref{picard} provides infinite towers of lattices in $\SU(2,1)$ with Property (FA) but not Property (T). The manifold assumption of Theorem \ref{GS} restricts us to the torsion free setting, however with Selberg's lemma we also have the following corollary.

\begin{cor}\label{secondfa}
Every congruence arithmetic lattice $\Gam<\SU(2,1)$ of second type has Property (FA).
\end{cor}

\begin{proof}
If $\Gam$ is torsion free, this is a direct application of Theorem \ref{picard}. If $\Gam$ has torsion, it suffices to show that $\Gam$ has a finite index normal subgroup with Property (FA) -- see Proposition \ref{normal}. Selberg's lemma implies that $\Gam$ has a finite index torsion-free normal subgroup $\Gam^\prime$, and $\Gam^\prime$ is a congruence subgroup by construction. Theorem \ref{picard} implies that $\Gam^\prime$ has Property (FA), so $\Gam$ must have Property (FA).
\end{proof}

One setting to which Theorem \ref{picard} applies is fake projective planes -- compact algebraic surfaces with the same betti numbers as $\C\mathbb{P}^2$. It follows from Yau's solution to the Calabai conjecture that all fake projective places are complex hyperbolic surfaces and that there are only finitely many up to homeomorphism (see \cite{Klin}). The first example was constructed by Mumford \cite{Mum}, and his construction implies that the corresponding lattice in $\PU(2,1)$ is of the second type arising from a congruence subgroup. In fact, it is proven independently in \cite{Klin,Yeung} that all fake projective planes are arithmetic, and more recently Prasad and Yeung \cite{PY} classified fake projective planes using arithmetic techniques. However, a fake projective plane need not be congruence of second type, so Theorem \ref{picard} may not apply. Nonetheless, using the archimedean superrigidity of fake projective planes \cite{Klin}, we will prove:

\begin{thm}\label{fakecp2}
The fundamental group of any fake projective plane does not split as a nontrivial free product with amalgamation. In particular, it has Property (FA).
\end{thm}

One way to think of fake projective planes is as complex hyperbolic cousins to rational homology 3-spheres. In contrast to Theorem \ref{fakecp2}, F. Calageri and N. Dunfield \cite{CD} construct, assuming certain conjectures in number theory (which are removed in \cite{BE}), an infinite tower of arithmetic hyperbolic rational homology 3-spheres $M_n$ such that the injectivity radius grows arbitrarily large as $n\to\infty$. However, unlike fake projective planes, $\pi_1(M_n)$ is Haken and splits as a nontrivial free product with amalgamation for all $n$ -- see \cite[Sec. 2.14]{CD}. In particular, the question of whether there are non-Haken hyperbolic 3-manifolds with arbitrarily large injectivity radius remains open. This analogy, as opposed to Theorem \ref{zeta3}, highlights one of the many differences between the real and complex hyperbolic worlds.

\section{Preliminaries}\label{prelims}

Here, we collect the definitions and facts necessary for later sections.

\subsection{Complex hyperbolic space}\label{ch2}

We briefly recall the construction of the complex hyperbolic plane $\Hy_\C^2$ and its finite volume quotients. Consider the Hermitian form on $\C^3$ of signature $(2,1)$ given by
\begin{equation}J=\begin{mat}0 & 0 & 1\\ 0 & 1 & 0\\ 1 & 0 & 0\end{mat},\end{equation}
which in coordinates is
\begin{equation}\langle z,w\rangle= z_1\conj{w}_3+z_2\conj{w}_2+z_3\conj{w}_1,\end{equation}
and let $N_-$ denote the collection of $z\in\C^3$ such that $\langle z,z\rangle<0$. Then, $\Hy_\C^2$ is the projectivization of $N_-$, which can be canonically identified with the open unit ball in $\C^2$ with the Bergman metric. It is clear from this construction that we obtain biholomorphic isometries of $\Hy_\C^2$ from the group
\begin{equation}\SU(J)=\{A\in\SL(3;\C)\ :\ A^*JA=J\},\end{equation}
where * denotes conjugate transposition. We will denote $\SU(J)$ by $\SU(2,1)$, though we should remark that this is somewhat nonstandard notation. This will be convenient for consistency with the notation of Falbel and Parker \cite{FP} that we require later.

Since we projectivize to obtain $\Hy_\C^2$, the group of biholomorphic isometries of $\Hy_\C^2$ is isomorphic to $\PU(2,1)$. Recall that $\SU(2,1)$ is a 3-fold cover of $\PU(2,1)$ by the subgroup generated by $\zeta_3I$, where $I$ is the identity matrix and $\zeta_3$ is a primitive cube root of unity. This allows us to identify the fundamental groups of complex hyperbolic surfaces and 2-orbifolds, finite volume quotients of $\Hy_\C^2$ by discrete groups of isometries, with lattices in $\SU(2,1)$.

\subsection{Arithmetic lattices in $\SU(2,1)$}

See \cite{McR} for a treatment of arithmetic lattices in $\SU(n,1)$; our presentation is heavily influenced by these notes. For $n=2$ there are two distinct constructions of arithmetic lattices, which we will call arithmetic lattices of the first and second type.

To construct arithmetic lattices of the first type, we start with a totally real number field $F$ of degree $n$ over $\Q$ and an imaginary quadratic extension $E/F$ with ring of integers $\mc{O}_E$ and Galois embeddings $\sig_1,\conj{\sig}_1,\dots,\sig_n,\conj{\sig}_n:E\lra\C$ such that $\sig_i|_F=\conj{\sig}_i|_F$. Then, choose an $E$-defined Hermitian matrix $H\in\GL(3;\C)$ such that $H$ has signature $(2,1)$ at $\sig_1$ and $\conj{\sig}_1$ and signature $(3,0)$ at $\sig_i$ and $\conj{\sig}_i$ for all $i\neq 1$. Finally, for an $\mc{O}_E$-order $\mc{O}$ we define
\begin{equation}\SU(H;\mc{O})=\{A\in\SL(3;\mc{O}):A^*HA=H\},\end{equation}
where $*$ again denotes conjugate transposition. Under equivalence of Hermitian forms we can associate $\SU(H,\mc{O})$ with a lattice in $\SU(2,1)$ under the isomorphism $\SU(H)\cong\SU(2,1)$. We call any lattice in $\SU(2,1)$ commensurable with some $\SU(H,\mc{O})$ an \emph{arithmetic lattice of the first type}.

\begin{ex}\label{picardbianchi}
Let $F=\Q$, $E=\Q(\sqrt{-d})$ for $d$ a square free natural number, and let $\mc{O}_d$ denote the ring of integers in $E$. We take $J$ as in Sec. \ref{ch2}, and then $\SU(J;\mc{O}_d)=\SU(2,1;\mc{O}_d)$ is an arithmetic lattice of the first type. Since these lattices contain unipotent elements, Godement's compactness criterion implies that $\SU(2,1;\mc{O}_d)$ is non-cocompact.
\end{ex}

Arithmetic lattices of the second type are constructed as follows. Again, choose a totally real number field $F$ and an imaginary quadratic extension $E/F$ with ring of integers $\mc{O}_E$. Also, choose a degree three Galois extension $L/E$ with $\Gal(L/E)=\langle\thet\rangle$ and let $K/F$ be the degree three totally real subfield of $L$. For an element $\al\in E$ such that
\begin{equation}N_{E/F}(\al)\in N_{K/F}(K^\times),\quad\quad \al\notin N_{L/E}(L^\times),\end{equation}
where $N_{k/k^\prime}$ denotes the field norm, we define the degree three cyclic algebra
\begin{equation}A=(L/E,\thet,\al)=\left\{\sum_{i=0}^2\beta_iX^i\ :\ X^3=\al, X\beta=\thet(\beta)X\ \textrm{for}\ \beta,\beta_i\in L\right\}.\end{equation}
A theorem of Wedderburn implies that this is a division algebra by our choice of $\al$. Also, by a theorem of Albert, our selection of $\al$ also ensures that $A$ admits an involution $\tau$ such that the restriction $\tau|_E$ from the natural inclusion $E\hookrightarrow A$ is complex conjugation. We call such an involution an \emph{involution of the second kind}, and we define a \emph{Hermitian element} of an algebra equipped with such an involution $\tau$ to be an element $h$ such that $\tau(h)=h$. Notice that this is precisely the usual notion of a Hermitian matrix when $h$ is a matrix and $\tau$ is conjugate transposition.

For a Hermitian element $h\in A$ and an $\mc{O}_E$-order $\mc{O}$ of $A$, we define
\begin{equation}\SU(h,\mc{O})=\{x\in\mc{O}: \tau(x)hx=h\}.\end{equation}
Then $A\otimes_E\C\cong\M(3,\C)$, and choosing $h$ so that it has signature $(2,1)$ under this tensor product we obtain an isomorphism $\SU(h)\cong\SU(2,1)$. Thus we can identify $\SU(h,\mc{O})$ with a lattice in $\SU(2,1)$, and we call any lattice commensurable with some $\SU(h,\mc{O})$ an \emph{arithmetic lattice of the second type}. Since $A$ is a division algebra, Godement's compactness criterion implies that all such lattices are cocompact.

\begin{ex}[Mumford's Fake $\C\mathbb{P}^2$ \cite{Mum}]
We will not construct Mumford's example; we only give the arithmetic construction commensurable with it (see \cite{PY,McR}). However, Mumford's construction also implies that the corresponding lattice in $\SU(2,1)$ arises from a congruence subgroup.

For $\zeta_7$ a primitive $7^{th}$ root of unity, $F=\Q$, $E=\Q(\sqrt{-7})$, and $L=\Q(\zeta_7)$, let $\lam=(-1+\sqrt{-7})/2$, $\al=\lam/\conj{\lam}$, and $\thet$ be the generator of $\Gal(L/E)$, which is given by $\zeta_7\mapsto\zeta_7^2$. Then, $A=(L/E,\thet,\al)$ has the involution of the second kind given explicitly by
\begin{equation}\tau(X)=\conj{\al}X^2,\quad\quad\tau(\beta)=\conj{\beta}\quad\textrm{for}\ \beta\in E.\end{equation}
Finally, define the Hermitian form $h$ and $\mc{O}_E$-order $\mc{O}$ in $A$ given respectively by
\begin{equation}h=\conj{\lam}X^2-\conj{\lam}X+(\lam-\conj{\lam})\end{equation}
\begin{equation}\mc{O}=\mc{O}_L\oplus\conj{\lam}X\mc{O}_L\oplus\conj{\lam}X^2\mc{O}_L.\end{equation}
Then, $\Gam=\SU(h,\mc{O})$ is an arithmetic lattice of the second type commensurable with Mumford's fake $\C\mathbb{P}^2$.
\end{ex}

\subsection{Congruence subgroups}

As above, let $E/F$ be an imaginary quadratic extension of a totally real number field and $G$ the algebraic group determined by the $E$-defined Hermitian form $H$ or division algebra $D$, according to whether we are constructing first or second type lattices, respectively. Let $\A_f$ denote the finite ad\`{e}les of $F$, and for any open compact subgroup $K$ of $G(\A_f)$ set $\Gam_K=G(F)\cap K<G(\R)\cong\SU(2,1)$, where $F$ is embedded in the ad\`{e}les via the diagonal embedding.

Then, $\Gam_K$ is an arithmetic lattice in $\SU(2,1)$ of the appropriate type, and we call such a subgroup a \emph{congruence subgroup}. Given an ideal $\mf{a}$ of $\mc{O}_F$, the diagonal embedding of $\mf{a}$ into the ad\`{e}les provides us with an open compact subgroup of $G(\A_f)$. When $K$ corresponds to such a subgroup, it is the kernel of reduction in $G(\mc{O}_F)$ modulo the ideal $\mf{a}$, and we call such a group a \emph{principal congruence subgroup}. It is not immediately clear from this definition, but congruence subgroups are precisely those containing some principal congruence subgroup.

\begin{rem}
The use of $F$ as opposed to $E$ is often a source of confusion. We use $F$ because we are concerned with real analytic structure arising from $\SU(2,1)$ as opposed to structure arising from $\SU(2,1)^\C=\SL(3;\C)$.
\end{rem}

\subsection{Property (FA)}

If $\mc{T}$ is a tree with an action by a group $G$, we denote by $\mc{T}^G$ the subtree of fixed points of the $G$-action. We say that $G$ has \emph{Property (FA)} if $\mc{T}^G\neq\emptyset$ for every tree $\mc{T}$ on which $G$ acts without inversions. The following theorem, which appears as \cite[Theorem 15]{Se}, is the fundamental theorem in the study of Property (FA).

\begin{thm}\label{bigfa}
A group $G$ has Property (FA) if and only if
\begin{enumerate}
\item $G$ is finitely generated.

\item $G$ does not split as a nontrivial free product with amalgamation.

\item $G$ does not admit a homomorphism onto $\Z$.
\end{enumerate}
\end{thm}

The following two propositions will be crucial in the proof of Theorem \ref{zeta3}. They appear as Example 6.3.3 on p. 60 and Proposition 26 on p. 64 of \cite{Se}, respectively, but we include their proofs for completeness.

\begin{prop}\label{normal}
Suppose $G$ is a finitely presented group and $N\trianglelefteq G$ a normal subgroup such that $N$ and $G/N$ have Property (FA), then $G$ also has Property (FA).
\end{prop}

\begin{proof}
Suppose $G$ acts on the tree $\mc{T}$ without inversions. Then $N$ acts on $\mc{T}$ and $\mc{T}^N\neq\emptyset$, and there is an induced action of $G/N$ on the tree $\mc{T}^N$ with $\mc{T}^\prime=(\mc{T}^N)^{G/N}\neq\emptyset$. $G$ fixes $\mc{T}^\prime$, so $G$ has Property (FA).
\end{proof}

\begin{prop}\label{subgroups}
Suppose $G$ is a group with subgroups $A=\langle a_i\rangle$ and $B=\langle b_j\rangle$ with $G=\langle A,B\rangle$ and that $G$ acts on a tree $\mc{T}$. If $\mc{T}^A,\mc{T}^B\neq\emptyset$ and every $a_ib_j$ has a fixed point on $\mc{T}$, then $\mc{T}^G\neq\emptyset$.
\end{prop}

\begin{proof}
Suppose $G$ acts without inversions on a tree $\mc{T}$. Since $G$ is generated by the subgroups $A$ and $B$, it follows that $\mc{T}^G=\mc{T}^A\cap\mc{T}^B$, so suppose that this intersection is trivial, i.e. $\mc{T}^A$ and $\mc{T}^B$ are disjoint inside $\mc{T}$. Then, we can find a nontrivial geodesic $\gam$ joining $\mc{T}^A$ and $\mc{T}^B$ with endpoints $x\in\mc{T}^A$ and $y\in\mc{T}^B$. If $x^\prime$ is the vertex of $\mc{T}$ on $\gam$ distance one from $x$, some generator $a_i$ of $A$ does not fix $x^\prime$, which implies that following $\gam$ from $y$ to $x$, then following $a_i\gam$ from $a_ix=x$ to $a_iy\neq y$ is a geodesic with midpoint $\gam\cap a_i\gam=x$.

Since $y\in\mc{T}^B$, we also have $a_iy=a_ib_jy$ for every $j$, so the geodesic from $y$ to $\mc{T}^{a_ib_j}\neq\emptyset$ to $a_ib_jy=a_iy$ must be $\gam$ followed by $a_i\gam$. Then, $a_ib_j$ fixes the midpoint of this geodesic, which is $x$, implying that $x\in\mc{T}^{a_ib_j}$ for all $j$. Therefore, $b_jx=a_i^{-1}x=x$ for all $j$, contradicting that $x\notin\mc{T}^B$.
\end{proof}

\section{The proof of Theorem \ref{zeta3}}\label{zeta3proof}

Let $\Gam_3$ denote the group $\SU(2,1;\mc{O}_3)$. Recall from Example \ref{picardbianchi} that this is the subgroup of $\SL(3;\mc{O}_3)$ preserving the Hermitian form
\begin{equation}J=\begin{mat} 0 & 0 & 1\\ 0 & 1 & 0\\ 1 & 0 & 0\end{mat}\end{equation}
of signature (2,1). Also, let $\mc{D}(\mc{O}_3)$ denote the diagonal subgroup of $\Gam_3$ and $\mc{N}(\mc{O}_3)$ denote the subgroup of strictly upper triangular matrices, which is a lattice in the 3-dimensional Heisenberg group. The Borel subgroup of upper triangular matrices is
\begin{equation}\mc{B}(\mc{O}_3)=\mc{N}(\mc{O}_3)\rtimes\mc{D}(\mc{O}_3),\end{equation}
and the Borel subgroup of $\PU(2,1;\mc{O}_3)$ is the projectivization of the Borel subgroup in $\SU(2,1;\mc{O}_3)$. That the Borel subgroup is preserved under projectivization follows from its characterization as the stabilizer of a point in the ideal boundary of each respective symmetric space. Similar to Serre's proof that $\PSL(2;\mc{O}_3)$ has Property (FA), we will make use of a particular presentation of $\PU(2,1;\mc{O}_3)$, which is stated as Theorem 5.9 of \cite{FP}.

\begin{thm}\label{falbelparker}
$\PU(2,1;\mc{O}_3)$ has a presentation
\begin{equation}\langle R,P,QP^{-1}\ :\ R^2,(QP^{-1})^6,(RP)^3,[R,QP^{-1}],P^3Q^{-2}\rangle,\end{equation}
where $\langle P,QP^{-1}\rangle$ generates the Borel subgroup.
\end{thm}

\begin{rem}
It is not stated explicitly in \cite{FP} that $\langle Q,P\rangle$ generates the Borel subgroup, but this follows immediately from their Proposition 3.2, where they prove that this group is the stabilizer of infinity for the action on the boundary of complex hyperbolic space, considered as Heisenberg space with the point at infinity.
\end{rem}

\begin{proof}[Proof of Theorem \ref{zeta3}]
First, we claim that the Borel subgroup, $\mc{B}(\mc{O}_3)$, has Property (FA). It follows immediately from the presentation that the Borel subgroup has finite abelianization, so it cannot map onto $\Z$. Indeed, the abelianization has $P^3=Q^2$ and $P^6=Q^6$, implying that $Q^4=Q^6$, so $Q^2=1$, which implies that $P^3=1$.

To show that it cannot split as a nontrivial free product with amalgamation, Proposition 3.1 in \cite{FP} shows that the Borel subgroup fits into a short exact sequence
\begin{equation}1\lra\Z\lra\mc{B}(\mc{O}_3)\lra\Del(2,3,6)\lra 1.\end{equation}
It follows from Proposition \ref{subgroups} that $\Del(2,3,6)$ has Property (FA), so it cannot split as a free product with amalgamation. Since the $\Z$ factor is central in $\mc{B}(\mc{O}_3)$, if $\mc{B}(\mc{O}_3)$ splits as a free product with amalgamation then the $\Z$ subgroup must be contained in the amalgamating subgroup. However, this implies that the short exact sequence induces a nontrivial free product with amalgamation for $\Del(2,3,6)$, which is a contradiction.

Now, we apply Proposition \ref{subgroups} to $\PU(2,1;\mc{O}_3)=\langle R,\mc{B}(\mc{O}_3)\rangle=\langle A,B\rangle$, where $\langle R\rangle\cong\Z/2\Z$ has Property (FA) since it is a finite group. In other words, given an action of $\PU(2,1;\mc{O}_3)$ on a tree $\mc{T}$, we know that $\mc{T}^A,\mc{T}^B\neq\emptyset$, so we need only prove that the products $RP$ and $R(QP^{-1})$ have fixed points on $\mc{T}$. This follows from the above presentation, since finite order elements necessarily have fixed points on $\mc{T}$ and
\begin{equation}(RP)^3=1,\quad\quad(RQP^{-1})^6=R^6(QP^{-1})^6=1.\end{equation}
Therefore, $\PU(2,1;\mc{O}_3)$ has Property (FA).

Finally, to show that $\SU(2,1;\mc{O}_3)$ has Property (FA), we apply Proposition \ref{normal} to the short exact sequence
\begin{equation}1\lra\Z/3\Z\lra\SU(2,1;\mc{O}_3)\lra\PU(2,1;\mc{O}_3)\lra 1\end{equation}
and the proof is complete.
\end{proof}

\section{The proofs of Theorems \ref{picard} and \ref{fakecp2}}\label{secondproof}

In order to prove Theorems \ref{picard} and \ref{fakecp2}, we will need some additional results about the K\"{a}hler structure of compact complex hyperbolic surfaces. Recall that a Riemannian manifold $(X,g)$ is a \emph{K\"{a}hler manifold} if it admits an integrable almost complex structure $J\in\End(TX)$ with $J^2=-\Id$ such that the form $\om(X,Y)=g(JX,Y)$ is closed. We call a group $\Gam$ a \emph{K\"{a}hler group} if it is the fundamental group of a compact K\"{a}hler manifold. In particular, all cocompact lattices in $\SU(n,1)$ are K\"{a}hler groups, as they give rise to complex projective varieties. The following striking theorem (\cite{GS} Theorem 9.1) connects the geometry of a K\"{a}hler manifold with the structure of its fundamental group.

\begin{thm}\label{GS}
Let $X$ be a compact K\"{a}hler manifold with fundamental group $\Gam=\Gam_1*_\Del\Gam_2$ where $\Del$ is of index at least 2 in $\Gam_1$ and of index at least 3 in $\Gam_2$, where either index is allowed to be infinite. Then $X$ maps holomorphically onto a compact Riemann surface.
\end{thm}

Let $H^{1,1}(X)$ denote the collection of 2-forms on a complex manifold $X$ that split into holomorphic and antiholomorphic part. The Picard number of $X$ to defined to be the rank of $H^{1,1}(X)\cap H^2(X,\Q)$. We say that a torsion-free lattice $\Gam<\SU(2,1)$ has Picard number one if the corresponding quotient manifold $\Hy^2_\C/\Gam$ has Picard number one. We will make use of the following lemma, due to Yeung \cite{Yeung}, whose proof we include for completeness.

\begin{lem}\label{Yeung}
If $X$ is an algebraic surface with Picard number one then $X$ admits no nontrivial holomorphic map onto a compact Riemann surface.
\end{lem}

\begin{proof}
Let $f:X\lra\Sig$ be a nontrivial holomorphic map from $X$ to a compact Riemann surface $\Sig$. Then, the fundamental class $[\Sig]$ pulls back to a non-torsion element $\sig\in H^{1,1}\cap H^2(X;\Z)$. Since $X$ has Picard number one, this is a nonzero multiple of the generator $\thet$ of $H^{1,1}\cap H^2(X;\Z)$, which implies that the push-forward of $\thet$ is a nontrivial cycle. Generic fibers of $f$ are one-dimensional over $\C$, so if $\al$ is the cohomology class representing a generic fiber it is also a nonzero multiple of $\thet$. Since $\thet$ has a nontrivial push-forward, $\al$ must also have a nonzero push-forward. However, generic fibers necessarily have trivial push-forward, which is a contradiction.
\end{proof}

We will also need the following elementary lemma.

\begin{lem}\label{index2}
Any normal infinite dihedral subgroup of $D_\infty$ has index at most 2.
\end{lem}

\begin{proof}
Recall that $D_\infty\cong\Z\rtimes(\Z/2\Z)$, so the infinite dihedral subgroups of $D_\infty$ are precisely $(\ell\Z)\rtimes(\Z/2\Z)$ for some $\ell\geq 1$. Consider $D_\infty$ acting on the real line by translations and negation, so that $(p,\ep)(x)=p+\ep x$ for $\ep=\pm 1$. Then, we have
\begin{equation}(p,-1)(p,-1)(x)=(p,-1)(p-x)=p-p+x=x,\end{equation}
for all $x\in\R$, so $(p,-1)=(p,-1)^{-1}$. Now, we conjugate the element $(\ell,-1)$ by $(p,-1)$ and see that
\begin{equation}(p,-1)(\ell,-1)(p,-1)(x)=2p-\ell-x=(2p-\ell,-1)(x).\end{equation}
This is in $(\ell\Z)\rtimes(\Z/2\Z)$ for all $p$ if and only if $\ell=1,2$. Since $(2\Z)\rtimes(\Z/2\Z)$ has index 2 in $D_\infty$, this proves the lemma.
\end{proof}

The following is Theorem 15.3.1 of \cite{Ro}.

\begin{thm}\label{Rog}
Let $\Gam$ be an arithmetic lattice in $\SU(2,1)$ of the second type arising from congruence and $X=\Hy^2_\C/\Gam$. Then
\begin{equation}H^1(X,\Q)=H^1(\Gam,\Q)=0.\end{equation}
In particular, $\rank(H_1(X,\Q))=b_1(X)=0$.
\end{thm}

The following theorem is often credited to \cite{Ro}, but the book contains no mention of such a result. In fact, for our lattices \cite{Ro} tells us about the cohomology in every dimension except 2. J. Rogawski kindly provided a copy of the correct reference, namely Theorem 3 of \cite{BR}.

\begin{thm}\label{BR}
If $\Gam$ is a congruence arithmetic lattice in $\SU(2,1)$ of the second type, then $\Gam$ has Picard number one.
\end{thm}

For the proof of Theorem \ref{picard}, we will need an extension of Theorem \ref{Rog} to the infinite dihedral group.

\begin{prop}\label{nodihedral}
If $\Gam<\SU(2,1)$ is a congruence arithmetic lattice of second type, then $\Gam$ admits no homomorphism onto the infinite dihedral group.
\end{prop}

\begin{proof}
Suppose that $\Gam_K$ is a congruence subgroup corresponding to the open compact subgroup $K<G(\A_f)$, and let $\mc{K}$ denote the collection of all open compact subgroups $K^\prime<G(\A_f)$ so that $\Gam_{K^\prime}$ is a principal congruence subgroup. Then, $K\cap K^\prime<G(\A_f)$ is an open compact subgroup for any $K^\prime\in\mc{K}$, and so $\Gam_{K\cap K^\prime}\triangleleft\Gam_K$ is a congruence subgroup. Our assumption that $\Gam_{K^\prime}$ is principal congruence assures that all subgroups that we consider are normal in $\Gam_K$.

Now, suppose that $\rho:\Gam_K\lra D_\infty$ is a surjective homomorphism with kernel $\Del$. Then, $\Gam_{K\cap K^\prime}\Del\triangleleft\Gam_K$ is a congruence subgroup, since it contains the congruence subgroup $\Gam_{K\cap K^\prime}$, with
\begin{equation}\Gam_{K\cap K^\prime}\Del/\Del\triangleleft\Gam_K/\Del\cong D_\infty.\end{equation}
Since the finite index subgroups of $D_\infty$ are isomorphic to either $\Z$ or $D_\infty$, and since $b_1(\Gam_{K\cap K^\prime})=0$ by Theorem \ref{Rog}, $\Gam_{K\cap K^\prime}\Del/\Del$ must be isomorphic to $D_\infty$ sitting normally inside $\Gam_K/\Del\cong D_\infty$. However, the only normal subgroups of $D_\infty$ that are isomorphic to $D_\infty$ have index at most 2 by Lemma \ref{index2}. Thus, to contradict the existence of $\rho$ it suffices to show that there exists some $K^\prime\in\mc{K}$ so that
\begin{equation}[\Gam_K/\Del:\Gam_{K\cap K^\prime}\Del/\Del]=[\Gam_K:\Gam_{K\cap K^\prime}\Del]>2.\end{equation}
To prove this, we first note that
\begin{equation}\bigcap_{K^\prime\in\mc{K}}\Gam_{K^\prime}=\{1\},\end{equation}
essentially since elements of $G(\mc{O}_F)$ are divisible by only finitely many primes. This implies that
\begin{equation}\bigcap_{K^\prime\in\mc{K}}\Gam_{K^\prime\cap K}\Del=\Del,\end{equation}
so that
\begin{equation}[\Gam_K:\bigcap_{K\cap K^\prime}\Gam_{K\cap K^\prime}\Del]=[\Gam_K:\Del]=\infty.\end{equation}

If $[\Gam_K:\Gam_{K\cap K^\prime}\Del]\leq 2$ for all $K^\prime\in\mc{K}$, it follows that $\Gam_{K\cap K^\prime}\Del$ necessarily lies in a finite list of finite index subgroups of $\Gam_K$ (possibly containing $\Gam_K$) arising from $\Hom(D_\infty,\Z/2\Z)$. However, this means that $\bigcap\Gam_{K\cap K^\prime}\Del$ is the intersection of finitely many finite index subgroups, which is necessarily of finite index in $\Gam_K$. This contradiction completes the proof of the proposition.
\end{proof}

Similarly, we will need to rule out representation onto the infinite dihedral group for the fundamental groups of fake projective planes.

\begin{prop}\label{nodihedral2}
The fundamental group of a fake projective plane cannot surject $D_\infty$.
\end{prop}

\begin{proof}
As mentioned in the introduction, this will follow from the archimedean superrigidity of these lattices, which is due to Klingler \cite{Klin}.

\begin{thm}\label{Klin}
If $\Gam$ is the fundamental group of a fake projective plane, any homomorphism $\rho:\Gam\lra\PGL(3;\C)$ has compact Zariski closure or extends to a totally geodesic embedding of $\PU(2,1)$ into $\PGL(3;\C)$.
\end{thm}

Now, we construct a family of faithful representations of $D_\infty$ into $\SL(3;\C)$ that factor through the inclusion $\GL(2;\C)\lra\SL(3;\C)$ given by
\begin{equation}A\mapsto\begin{mat}A & 0\\ 0 & \frac{1}{\det{A}}\end{mat}\end{equation}
and that have eigenvalues off the unit circle $S^1$. To complete the proof with such a representation $\rho$, let $\conj{\rho}:\Gam\lra\SL(3;\C)$ be the representation obtained by composing the natural surjection $\Gam\lra D_\infty$ with $\rho$. Since we can choose $\conj{\rho}$ to have arbitrarily large eignevalues, it follows that $\conj{\rho}(\Gam)$ cannot have compact Zariski closure. It also follows that it does not arise from a totally geodesic embedding $\SU(2,1)\lra\SL(3;\C)$, since this would produce a totally geodesic embedding of $\SU(2,1)$ in $\GL(2;\C)$, which is impossible. We then projectivize this representation to obtain the same results in $\PGL(3;\C)$. This contradicts Theorem \ref{Klin} and completes the proof.

To construct the representation of $D_\infty$, present $D_\infty$ as $\langle r,s\ :\ s^2, srsr\rangle$ and consider the matrices in $\GL(2;\C)$ given by
\begin{equation}R=\begin{mat}\al & \beta\\ 2\im(\al)i & \conj{\al}\end{mat}\quad\quad S=\begin{mat}1 & -1\\ 0 & -1\end{mat}.\end{equation}
A pair of calculations shows that $S^2=I$ and
\begin{equation}SRS^{-1}=\det(R)R^{-1},\end{equation}
so if $\det(R)=\al\conj{\al}-2\im(\al)\beta i=1$, we obtain a representation $\rho$ of $D_\infty$ into $\GL(2;\C)$. Furthermore, if the eigenvalues of $R$ lie off the unit circle, it follows that this representation will be faithful and that the image cannot lie in a conjugate of the unitary group, and thus does not have compact Zariski closure. Since the equation for the eigenvalues of $R$ is
\begin{equation}1-2\real(\al)\lam+\lam^2=0,\end{equation}
we can obtain any nonzero eigenvalue $\lam_0$ we like by selecting
\begin{equation}\real(\al)=\frac{1+\lam_0^2}{2\lam_0},\end{equation}
as long as this number lies in $\R$. As $\beta$ and $\im(\al)$ do not factor into this equation, we still have the necessary freedom to assure that $\det(R)=1$. This provides us with the representation $\rho$ required above.
\end{proof}

Now, we are ready to prove Theorems \ref{picard} and \ref{fakecp2}.

\begin{proof}[Proof of Theorems \ref{picard} and \ref{fakecp2}]
Let $\Gam<\SU(2,1)$ be a torsion-free cocompact lattice satisfying the hypotheses of either theorem, and suppose that $\Gam$ splits as a nontrivial free product with amalgamation $\Gam_1*_\Del\Gam_2$. Theorem \ref{GS} combined with Lemma \ref{Yeung} allows us to assume that $[\Gam_i:\Del]=2$ for $i=1,2$. To see this, notice that if $\Del$ has index at least 3 (possibly infinite) in either of the $\Gam_i$, Theorem \ref{GS} gives a holomorphic map onto a compact Riemann surface, which contradicts Lemma \ref{Yeung}. Also, notice that $\Gam=\Gam_1*_\Del\Gam_2$ with $[\Gam_i:\Del]=2$ for $i=1,2$ if and only if $\Gam$ surjects the infinite dihedral group $D_\infty$. This is ruled out by Proposition \ref{nodihedral} for the congruence lattices and by Proposition \ref{nodihedral2} for fake projective planes.
\end{proof}

\section*{Acknowledgments}
First, I would like to thank my advisor, Alan Reid, for raising these questions and for many invaluable conversations and careful readings of this paper. Also, thanks to Lior Silberman for suggestions, Jon Rogawski for providing \cite{BR}, Ben McReynolds for many helpful discussions, Domingo Toledo for helpful correspondence and for pointing out that an early version of this paper used a result of Reznikov whose proof contains an error, and the referee for a very attentive reading. I also wish to thank the Centre Interfacultaire Bernoulli of EPF Lausanne for their hospitality during the preparation of part of this work.


\end{document}